\documentclass{amsart}
\usepackage{graphicx}
\usepackage{amssymb,epic,eepic}

\newtheorem{thm}{Theorem}[section]
\newtheorem{lem}[thm]{Lemma}

\newtheorem{cor}[thm]{Corollary}

\newtheorem{rem}[thm]{Remark}
\newtheorem{quest}[thm]{Question}
\newtheorem{conj}[thm]{Conjecture}

\newcommand{\Si}{{\Sigma}}
\newcommand\tM{{\widetilde M}}
\newcommand\tS{\widetilde \Sigma}
\newcommand\bS{\Sigma_b}
\newcommand\ra{\rightarrow}

\newcommand\mapright[1]{\smash{\mathop{\longrightarrow}\limits^{#1}}}
\newcommand\abs[1]{|#1|}

\newcommand\R{{\mathbb R}}
\newcommand\Z{{\mathbb Z}}

\newcommand\C{{\mathbb C}}

\DeclareMathOperator{\Tr}{Tr}

\def\Id{{\mathop{\fam0 Id}\nolimits}}
\def\SL{{\mathop{\fam0 SL}\nolimits}}
\def\GL{{\mathop{\fam0 GL}\nolimits}}
\def\PSL{{\mathop{\fam0 PSL}\nolimits}}

\newcommand\ta{
\begin{picture}(10,10)
\path(0,0)(10,0)(10,10)(0,10)(0,0)
\path(5,0)(5,10)
\path(0,5)(10,5)
\end{picture}}

\newcommand\tatableauxf{\begin{picture}(20,20)
\path(0,0)(20,0)(20,20)(0,20)(0,0)
\path(10,0)(10,20)
\path(0,10)(20,10)
\put(2,2){$1$}\put(12,2){$2$}\put(2,12){$3$}\put(12,12){$4$}
\end{picture}}

\newcommand\tatableauxs{\begin{picture}(20,20)
\path(0,0)(20,0)(20,20)(0,20)(0,0)
\path(10,0)(10,20)
\path(0,10)(20,10)
\put(2,2){$1$}\put(12,2){$3$}\put(2,12){$2$}\put(12,12){$4$}
\end{picture}}

\begin{document}

\title[TQFT and the Nielsen-Thurston
  classification of $M(0,4)$]{Topological Quantum Field Theory and the Nielsen-Thurston
  classification of $M(0,4)$}
\author{J{\o}rgen Ellegaard Andersen}
\address{Department of Mathematics\\
        University of Aarhus\\
        DK-8000 Aarhus C, Denmark}
\email{andersen@imf.au.dk}

\author{ Gregor Masbaum}
\address{Institut de Math{\'e}matiques de Jussieu (UMR 7586 du CNRS)\\
Universit{\'e} Paris 7 (Denis Diderot) \\
Case 7012\\
2, place Jussieu\\
75251 Paris Cedex 05\\
FRANCE }
\email{masbaum@math.jussieu.fr}

\author{Kenji Ueno}
\address{Department of Mathematics\\
        Faculty of Science, Kyoto University\\
        Kyoto, 606-01 Japan}
\email{ueno@math.kyoto-u.ac.jp}

\maketitle

\begin{center}{\em to appear in:} Mathematical Proceedings of the
    Cambridge Philosophical Society \end{center}

\begin{abstract}
We show that the Nielsen-Thurston classification of mapping classes of
the sphere with four marked points is determined by the quantum
$SU(n)$ representations, for any fixed $n\geq 2$.  In the
Pseudo-Anosov case we also show that the stretching factor is a limit
of eigenvalues of (non-unitary) $SU(2)$-TQFT representation matrices.  It follows
that at big enough levels, Pseudo-Anosov mapping classes are represented by matrices of
infinite order.

\end{abstract}

\tableofcontents

\section{Introduction} \label{sec.intro}

Quantum invariants of $3$-manifolds and the associated
`Topological Quantum Field Theories' (TQFT) give rise to
finite-dimensional representations of  mapping class groups of
surfaces. It is
an obvious quest to seek a relationship
between these representations and the Nielsen-Thurston
theory of mapping classes.

The representations can be constructed using a number of different
techniques. (See also \cite{SMF} for a survey.) On the combinatorial side one has the constructions of Reshetikhin
and Turaev \cite{RT1} and \cite{RT2} based on the theory of
quantum groups, and more generally by Turaev based on the theory of
modular categories. Using skein theory there are the constructions
of Blanchet, Habegger, Masbaum and Vogel
\cite{BHMV2} based on the Kaufmann bracket and \cite{Bl} using the
HOMFLY-polynomial.

In the geometric approaches the representations come from
(projectively) flat bundles over Teichm{\"u}ller spaces.
The construction of these bundles with their projectively flat
connections was first given from the point of view of conformal
field theory by Tsuchiya, Ueno and Yamada \cite{TUY}, \cite{Ue} and
\cite{Ue2} using the theory of representations of affine Lie
algebras.
In \cite{AU1}, \cite{AU2}
and \cite{AU3} (in preparation), the first and last author are
currently working
to establish that these representations are
isomorphic to the ones constructed by the skein methods of
\cite{BHMV2}. Using the method of
geometric quantization of moduli spaces of flat connections
one gets yet another
construction, which is due to Axelrod, Della Pietra and Witten \cite{ADW}
and Hitchin \cite{H}.\footnote{That these bundles are isomorphic to
the ones constructed from conformal field theory is due to Laszlo
\cite{La1}.}  In fact Teichm{\"u}ller space naturally
parametrizes K{\"a}hler structures on the moduli space of flat
connections and the fiber of the bundle over a point in
Teichm{\"u}ller space is the geometric quantization of the
moduli space with respect to the corresponding K{\"a}hler
structure.

Motivated by the above quest, the first author analyzed in \cite{A2} (and
in \cite{A1}) the behavior of the family of K{\"a}hler structures
on the moduli space near the Thurston boundary of Teichm{\"u}ller
space and found many new polarizations on the moduli
spaces in this way. The relevance of this being that the Thurston
boundary
of Teichm{\"u}ller space has the
following two fundamental properties, namely it compactifies
Teichm{\"u}ller spaces to a topological closed ball and the action
of the mapping class group extends to a continuous action on the
compactification, so therefore any
mapping class has a fixed point in Thurston's compactification of
Teichm{\"u}ller space. Points on the boundary are
represented by projective classes of measured foliations on the
surface and in fact by the Nielsen-Thurston classification, a mapping class
which is neither finite order nor reducible, is a Pseudo-Anosov
mapping class, that is, it preserves two transverse foliations.
These two foliations support unique (up to scale) transverse
measures and the Pseudo-Anosov mapping class scales one by the
so-called stretching factor and the other by one over the stretching
factor.

This leaves us with the following two very natural
questions (see also problem 8.11 formulated by the second
author in Ohtsuki's problem list \cite{Oh}):

\renewcommand{\thefootnote}{\fnsymbol{footnote}}

\begin{quest}\label{q1} \mbox{ }

(1) Can one determine the Nielsen-Thurston classification of mapping
classes using quantum representations? \footnote{{\em Note added in
    proof:} For an affirmative answer to Question 1.1.(1) see \cite{A4}.}

(2) If so, can one compute the stretching factor for Pseudo-Anosov mapping classes from the quantum representations?

\end{quest}

\addtocounter{footnote}{-1}
\renewcommand{\thefootnote}{\arabic{footnote}}

The main result of this paper is to give a positive answer to these
questions in the case of the mapping class group $M(0,4)$ of a surface
of genus zero with four marked points.

\begin{thm} For every $n\geq 2$, the Nielsen-Thurston classification
  of mapping classes in  $M(0,4)$
are determined by the quantum $SU(n)$ representations.
Moreover in the Pseudo-Anosov case, these representations also determine the stretching factor.
\end{thm}

This will be proved in Section~\ref{sec3}, using the geometric
approach to the representations. We also note an analog of
the asymptotic faithfulness result of \cite{A3} (see also \cite{FWW})
in this case
(Corollary~\ref{asf}). In Section~\ref{sec4},
we discuss the skein-theoretical approach (although for simplicity in the
$SU(2)$ case only) and obtain the following sharpening of our result:

\begin{cor}  A mapping class $\phi\in M(0,4)$ is
  Pseudo-Anosov if and only if its
$SU(2)$-TQFT representation matrix
  has infinite order for large enough level $k$.
Moreover, if $\phi$ is
  Pseudo-Anosov with stretching factor $\lambda(\phi)$, then for $k$
  big
enough the $SU(2)$-TQFT matrix
$\rho_k^{(S)}(\phi)$ has a unique eigenvalue $\lambda_k$ such that
$\vert\lambda_k\vert >1$, and \[ \lim_{k\ra \infty} \vert\lambda_k\vert
=  \lambda(\phi)~.\]
\end{cor}

Here, we must choose the roots of unity in the definition of the $SU(2)$-TQFT
appropriately (in the second statement only). We remark that for
this choice, the TQFT
is non-unitary (for $k$ big enough). See Section~\ref{sec4}
for details.
\vskip 8pt

\noindent{\em Acknowledgements.}

This work was done during visits of the last two authors
to Aarhus University and of the first two authors to Kyoto University.
We thank both institutions for their hospitality. We also thank
N. A'Campo for stimulating
discussions.

\section{A motivating example, train tracks, and open questions.}

Let  $\Si=\Si_{0,4}$ be a surface of genus zero with four marked points
$p_1, p_2,p_3, p_4$. We denote its mapping
class group $M(0,4)$ simply by $M$. Thus $M$ is the group
of path components of the group of orientation preserving self-homeomorphisms of $\Si$ which map the set of marked points to itself.

By blowing up each marked point to a boundary circle, we obtain a
compact surface $\bS$ with $4$ boundary components. We can also think
of $M$ as the isotopy classes of orientation preserving
self-homeomorphisms of $\bS$ where neither the homeomorphisms nor the
isotopies are required to be the identity on the boundary. This point of view is
the natural one for describing the measured
foliations associated to pseudo-Anosov mapping classes (see {\em
  e.g.} \cite{FLP}).

We have the following well known presentation of $M$ in terms of the standard
generators of the braid group $B_4$ (see Fig.~\ref{Fig1} below).

\begin{figure}[h]
  \begin{center}
    \includegraphics[scale=.4]{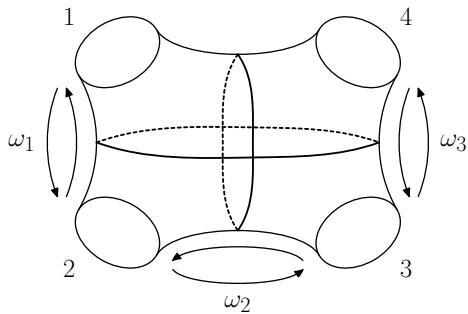}
    \caption{\small The braid generators $\omega_i$ and the Penner system.}
  \label{Fig1}
  \end{center}
\end{figure}

\begin{thm}\cite[Theorem 4.5]{B} \label{2.1}
The mapping class group $M=M(0,4)$ has a presentation with
generators $\omega_1, \omega_2, \omega_3$ and relations
\begin{eqnarray*}
\omega_1\omega_3 &=& \omega_3\omega_1\\
\omega_1 \omega_{2}\omega_{1} &=& \omega_{2} \omega_{1}\omega_{2}\\
\omega_{2} \omega_{3}\omega_{2} &=& \omega_{3} \omega_{2}\omega_{3}\\
\omega_1\omega_2\omega_3^2\omega_2\omega_1 &=& 1\\
(\omega_1\omega_2\omega_3)^4 &=& 1
\end{eqnarray*}
\end{thm}

\begin{figure}[h]
  \begin{center}
     \includegraphics[scale=.4]{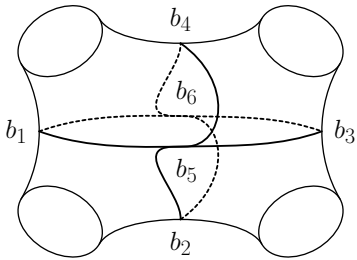}
    \caption{\small The train track $\tau$.}
    \label{Fig2}
  \end{center}
\end{figure}

\noindent{\bf Example.}
In \cite{P} Penner gave a method for
constructing Pseudo-Anosov mapping classes from certain systems of
curves.
We consider the system of two curves
illustrated in Fig.~\ref{Fig1} each of which cuts $\Si$ into two
pairs
of pants. By Penner's construction this system gives the train track $\tau$ on $\Si$ illustrated in Fig.~\ref{Fig2}.

This train track has $6$ branches and we label them $b_1$ to $b_6$ as indicated in Fig.~\ref{Fig2}.
A transverse measure $\mu$ on the train track $\tau$ is an assignment $\mu(b_i) = \mu_i$,
of non negative real numbers to the branches of $\tau$, which satisfies the switch conditions.
Hence we see that $\mu_5$ and $\mu_6$ are determined by $\mu_1,
\ldots, \mu_4$ and that they have to satisfy
\[\mu_1 + \mu_2 = \mu_3 +\mu_4.\]
The mapping classes $\omega_1^{-1}$, $\omega_2$ and $\omega_3^{-1}$ preserve $\tau$ and give the following incidence matrices
\[I(\omega_1^{-1}) = \left( \begin{array}{cccc} 1 & 0 & 0& 0 \\ 1 & 1 &0&0\\ 0& 0 & 1 & 0 \\ 1 & 0 & 0 & 1
  \end{array}\right)~,\  I(\omega_2) = \left( \begin{array}{cccc} 1& 1 & 0& 0 \\ 0 & 1 &0&0\\ 0& 1 & 1 & 0 \\ 0 & 0 & 0 & 1
  \end{array}\right)~,\  I(\omega_3^{-1}) = \left( \begin{array}{cccc} 1 & 0 & 0& 0 \\ 0 & 1 &1&0\\ 0& 0 & 1 & 0 \\ 0 & 0 & 1 & 1
  \end{array}\right). \]

Suppose now $\phi\in M$ is represented  by any positive word $w$ in $\omega_1^{-1},\omega_2$ and $\omega_3^{-1}$, which contains at least one $\omega_2$ and one $\omega_1^{-1}$ or one $\omega_3^{-1}$.
Then $\phi$ preserves the train track $\tau$ and its incidence matrix
$I(\phi)=I(w)$  is
Perron-Frobenius, so that all such $\phi\in M$ are Pseudo-Anosov. \footnote{In
  fact any Pseudo-Anosov mapping class in $M$ is conjugate to one of
  this form modulo the translation subgroup $N$ (see
  Section~\ref{secNT}). This follows from
Matsuoka's classification of Pseudo-Anosov conjugacy classes in the
braid group on 3 strands \cite{Mat}.}
In particular, $I(\phi)$ has a
unique maximal eigenvalue $\lambda >1$ which is the
stretching factor of $\phi$.

Now consider the representation matrices of $\phi$ in TQFT. As we will
prove, $\lambda$ is a limit of (absolute values of)
eigenvalues of such matrices. We can therefore ask:

\begin{quest} What is the relation between (any of) the constructions of
TQFT representations and  measured train tracks?
\end{quest}

At the time of this writing, we don't have a good answer to this
question. Our proofs exploit, instead,
the fact that the Nielsen-Thurston classification of $M(0,4)$ is in
some sense homological; for example, the stretching factor of $\phi$
can be computed from the action of $\phi$ on the homology of the
double cover of the sphere branched over the four marked points. (See Section~\ref{secNT}.)

While appropriate for $M(0,4)$, this reasoning does not
generalize to higher genus, as the stretching factor of a
Pseudo-Anosov mapping class in general cannot
be computed by homological means.
This leads us to mention the
following related questions which have puzzled us for some time.

\begin{quest} 1. Can one make sense of factorizing quantum representations
along a  measured train track?

2. Can one define the space of conformal blocks on a surface
equipped with a measured foliation as opposed to a
conformal structure?

3. Does a generic measured foliation determine a Lagrangian foliation
on the moduli space of flat connections?

4. What is the dynamics of these flat bundles over Teichm{\"u}ller space
with respect
to the flow on Teichm{\"u}ller space induced by a measured foliation?
\end{quest}

We end this section with the following
\begin{conj} A Pseudo-Anosov mapping class (on a surface of arbitrary
  genus) is represented in TQFT by
  matrices of infinite order, except for finitely many values of the
  level $k$.
\end{conj}

\section{The Nielsen-Thurston classification of mapping classes in
  $M(0,4)$.}\label{secNT}

In this section, we review the Nielsen-Thurston classification of mapping classes in
  $M= M(0,4)$ and remark that it is determined by the ``homology representation'' $h: M
  \ra \PSL_2(\Z)$.

Let $\pi : \tS \ra \Si$ be the double branched
cover of the sphere $\Si$ ramified at the four marked points. We can choose an identification of $\tS$ with the $2$-torus $T=\R^2/\Z^2$
such that the covering transformation is minus the identity acting on
the
torus. Recall that the mapping class group of the torus, $M(1,0)$, is
$\SL_2(\Z)$. We can describe the mapping class group $M=M(0,4)$ in a
similar fashion.

Given a matrix $A\in \SL_2(\Z)$
  and a vector $v\in (\frac 1 2 \Z)^2 \subset \R^2$, the
  transformation \[x \mapsto A x + v\] of $\R^2$ defines a
  diffeomorphism of $T=\R^2/\Z^2$ which commutes with $-\Id$. We
  denote by $\phi_{A,v}$ the induced diffeomorphism
  of $\Si= T/({-\Id})$.

\begin{thm}\label{2.2} All mapping classes in $M$ are represented by
  diffeomorphisms of the form $\phi_{A,v}$.
\end{thm}

\begin{proof} It suffices to check that the
  braid generators $\omega_i$ (see Theorem~\ref{2.1})  can be
  given in this way. Indeed, we can represent them as \[\omega_i=\phi_{A_i,v_i} \ \ \ \ \ \ \   (i=1,2,3)\]
  where \[A_1=A_3=\left( \begin{array}{cc} 1& 1\\ 0 & 1
  \end{array}\right),  \  A_2= \left( \begin{array}{cc} 1& 0\\ -1 & 1
  \end{array}\right),  \  v_1={{1/2}\choose 0}, \  v_2=
  {0 \choose {1/2}}, \  v_3={0\choose 0}~.\]
\end{proof}
\begin{rem}{\em In this description, the four ramification points
    $p_1, p_2,p_3,p_4$ are represented by the points in $\R^2$ with
coordinates ${{1/2}\choose 0},{0\choose 0},{0\choose {1/2}},
{{1/2}\choose {1/2}}$, respectively.}\end{rem}

Of course, $A$ and $v$ are not uniquely determined by $\phi_{A,v}$. Analyzing the
ambiguity one obtains the following
\begin{cor} \label{3.3} We have a short exact sequence
\[ 1\ra N \ra M\mapright{h} \ \PSL_2(\Z)\ra 1~,\] where $ h(\phi_{A,v})=
\pm A\in\PSL_2(\Z)$, and \[N \cong \Z/2\Z \times \Z/2\Z~.\]
\end{cor}

We call $h$ the {\em homology representation} of
$M$, because $\pm A$ is already determined by the action on the
homology of the torus $T$. We also call $N$ the {\em translation
  subgroup} of $M$.

\begin{proof}[Proof of Corollary \ref{3.3}.]  Let $M'$ be the group of diffeomorphisms (not
  considered up to isotopy) of the form
  $\phi_{A,v}$. It is easy to see that the corollary holds with $M'$
  in place of $M$. It remains to see that $M'=M$. In other words, we
  must show that two
  elements of $M'$ are isotopic if and only if they coincide. This can
  be checked by looking at the action on the homology
  of $T$ (recall this is given by $h$) and how the four ramification
  points are permuted.  \end{proof}

\begin{rem} {\em It is clear by elementary topological
considerations that any mapping class in $M$ lifts to a
self-homeomorphism of $T$, well-defined up to the covering transformation. As the mapping class group of the torus is
$\SL_2(\Z)$, this provides an alternative definition of $h : M\ra
\PSL_2(\Z)$ which does
not invoke
Theorem~\ref{2.2}.} \end{rem}

\begin{rem}{\em One can check that $M$ is a semi-direct product of $\PSL_2(\Z)$ with
$N$, the action of
$\PSL_2(\Z)$ on $N$ being the natural one {\em via} the homomorphism
$\PSL_2(\Z)\ra \PSL_2(\Z/2\Z)$. Birman \cite[section
 5.4]{B} discusses a similar semi-direct product decomposition
 starting from the presentation of $M$ given in Theorem~\ref{2.1}.
}\end{rem}

We recall that a matrix $A\in \SL_2(\Z)$ (or an element $\pm A\in \PSL_2(\Z)$)
is called {\em Anosov} if and only if
\[|\Tr(A)| > 2~.\] In that case,  the matrix $A$ has real eigenvalues
(w.l.o.g.) $\lambda >1$
and $\lambda^{-1}$, and for any choice of $v\in  (\frac 1 2 \Z)^2$,
the associated mapping class $\phi_{A,v}$ is
pseudo-Anosov with stretching factor $\lambda$.
 We briefly recall the
construction.  The action of $A$ on $T=\R^2/\Z^2$ preserves two
linear foliations
with irrational slope, which are also preserved by the translation by
$v$. The induced foliations on $\Si$ have
singularities at the four ramification points, but give rise to the right
kind of measured foliation on the blown up surface with boundary,
$\Si_b$.  (Compare
\cite[p. 217]{FLP}.\footnote{But notice that contrary to \cite{FLP} we
  allow
  permutations of the ramification points.} )

Conversely, if $\phi=\phi_{A,v}$ is  pseudo-Anosov, it is easy to see
that $h(\phi)=\pm A$ must be Anosov.
The following lemma is now almost obvious.

\begin{lem} \label{L1}
A mapping class $\phi\in M$ is finite order, reducible or Pseudo-Anosov if and only if $h(\phi)\in \PSL_2(\Z)$ respectively is finite order, reducible or
Anosov. Moreover, for a pseudo-Anosov mapping class $\phi$ we have that
 the stretching factor $\lambda(\phi)$ is $\abs{\lambda}$, where $\pm \lambda$  is
the (real) eigenvalue of the matrix $\pm h(\phi)$ such that
$\abs{\lambda}>1$.
\end{lem}
\begin{proof}
First observe that $\phi\in M$ is finite order if and only if
$h(\phi)$ is, since $N=\ker(h)$ is finite. The preceding
discussion shows that $\phi$ is pseudo-Anosov if and only if $h(\phi)$
is Anosov.
By exclusion according to the Nielsen-Thurston classification of diffeomorphisms (see e.g. \cite{FLP}), the lemma follows.
\end{proof}

\section{The quantum $SU(n)$ representations determine the homology representation}\label{sec3}

Choose a projective tangent vector at each of the marked points $p_i\in\Si$.
Let $\tM$
be the group of components of the group of orientation preserving self-diffeomorphisms
 of $\Si$ which map the set consisting of the four chosen projective tangent vectors to itself.
 Then $\tM$ is a quotient of the ribbon braid group
on four strands $RB_4$.
We have of course also a surjection $\tM\ra M$ whose kernel is generated by the
Dehn-twists around the four marked points. Let $\sigma_1, \sigma_2$ and $\sigma_3$ be the standard (i.e. no twists in the bands that make up the $\sigma_i$)
lifts to $RB_4$ of $\omega_1, \omega_2$ respectively
$\omega_3$. Pick a fifth point with a projective tangent vector in the complement of the four marked points on $\Si$.
Denote this point $\infty$ and let $\Si'$ be the resulting
surface with these five marked points. Then we see $RB_4$ as the subgroup of the ribbon mapping class group of $\Si'$ which fixes $\infty$.

Let $k\geq 2$ be an integer and consider the $SU(n)$-TQFT at level
k. The label set of this theory is the set of level $k$
integrable highest weight representations of the affine
Lie algebra associated to the Lie algebra of $SU(n)$, which is indexed
by
the
set
of Young diagrams with less than $n$ rows and less than or equal to $k$ columns. There is also an involution $\dagger$ of this set, which takes
the dual representation. - This means in particular
that if we label the marked points of $\Si'$ by such diagrams, we get a
finite dimensional vector space associated to this surface. If we choose the same label on the first four marked points,
then $RB_4$ acts on this vector space.
We now choose the label $\Box$ at four marked points of $\Si$ and
 $\ta^{\dagger}$ (the dual of $\ta$)  at $\infty$. For $SU(2)$
this diagram is not in the label set, but in this case it represents the trivial representation.
Let $V_{n,k}$ denote the resulting vector space and $\tilde\rho_{n,k}$ the resulting representation of $RB_4$.
For the construction of $V_{n,k}$ and $\tilde \rho_{n,k}$ in terms of
conformal blocks we refer to \cite{TUY}, \cite{Ue2} and \cite{K}.
For the skein theory construction we refer to \cite{BHMV2} and
\cite{B}. In this section we will focus on the geometric definition.

Let $q= e^{2 \pi i/(k+n)}$.

As described in \cite{K} (and \cite{TK} for the $SU(2)$-case) there is the path basis
for $V_{n,k}$ indexed by Young tableaux on $\ta$~. There are exactly two such tableaux, namely
\[\tatableauxf \ \ \  \mathrm{and} \ \ \  \tatableauxs ~.\]
Scale the vector corresponding to the first tableau by $B^{1/2}$ and the other vector by $B^{-1/2}$, where $B = \sqrt{q^3+q^2+q}$.
From the explicit formulae in \cite{K} one then gets in this rescaled path basis that
\[\tilde\rho_{n,k}(\sigma_1) =  \tilde\rho_{n,k}(\sigma_3) = q^{-\frac{n+1}{2n}}\left( \begin{array}{cc} q& 0\\ 0 & -1
  \end{array}\right)\]
\[ \tilde\rho_{n,k}(\sigma_2) = q^{-\frac{n+1}{2n}}\frac{1}{1+q}\left( \begin{array}{cc} -1& q^3+q^2+q\\ 1 & q^2
  \end{array}\right).\]

Moreover, each of the boundary parallel Dehn-twists (around each of the first four
points) acts by $q^{\frac{n+1}{2n}}$ times the identity.

We will see (by the computation below) that $\tilde\rho_{n,k}$ is
actually a projective
representation of our original mapping class group $M$. Here we use the following definition of a projective
representation of a group given by generators and relations: it is an
assignment of matrices to the generators such that the relations are
satisfied up to scalar multiples of the identity
matrix.

\begin{rem}{\em In the case $n=2$ we know {\em a priori} that
$\rho_{2,k}$ defines at least a projective representation of $\tM$,
and, hence, one of $M$ (since the kernel of the projection $\tM\ra M$ acts
 by
scalar matrices). This is because the label at $\infty$ is trivial for
$SU(2)$.
}\end{rem}

Let us now relate $\tilde\rho_{n,k}$ to the homology representation
$h$. First, we conjugate $\tilde\rho_{n,k}$ by \[C = \left( \begin{array}{cc} 1& q^2+q+1\\ 0 & q+1
  \end{array}\right)\]
to get upper and lower triangular matrices, as follows:
\[C\tilde\rho_{n,k}(\sigma_1)C^{-1} = C \tilde\rho_{n,k}(\sigma_3)C^{-1} = q^{-\frac{n+1}{2n}}\left( \begin{array}{cc} q& -(q^2+q+1)\\ 0 & -1
  \end{array}\right)
  \]
\[ C\tilde\rho_{n,k}(\sigma_2)C^{-1} = q^{-\frac{n+1}{2n}}\left( \begin{array}{cc} q& 0\\ 1 & -1
  \end{array}\right) .\]
 Motivated by this, we consider the following matrices $A,B$ with
 coefficients in the ring of Laurent polynomials $\Z[s,s^{-1}]$:
\[A = s^{-1}\left( \begin{array}{cc} s^2& -(s^4+s^2+1)\\ 0 & -1
  \end{array}\right) , \ \ \ \ \ \  B= s^{-1}\left( \begin{array}{cc} s^2& 0\\ 1 & -1
  \end{array}\right)~.\]

Then one can check that $ABA=BAB$ and $(AB)^3=1$, so that
\[\rho(\omega_1) =  \rho(\omega_3) = A, \ \ \ \  \rho(\omega_2) = B\]
defines a linear representation $\rho : M \ra \GL_2(\Z[s,s^{-1}])$. We now define
\[\rho_{n,k}(\omega_i)= q^{1/2n} C \tilde\rho_{n,k}(\sigma_i)
C^{-1}~,\] then, by construction, we have
\[ \rho_{n,k}(\omega_i) = \rho(\omega_i)|_{s=q^{1/2}}~.\]

Thus $\rho_{n,k}$ is a linear representation of $M$. It is also
clearly equivalent, as a projective representation, to the original
$\tilde \rho_{n,k}~.$
Moreover, since Laurent polynomials are
determined by their values at infinitely many points, the rescaled
TQFT representations $\rho_{n,k}$
determine $\rho$ uniquely and {\em vice versa}. (Note that it
 is enough to
consider them for a fixed $n$.) We therefore call
$\rho$ the universal quantum representation of $M$.

It remains to relate $\rho$ to the homology representation $h$. For
this, we evaluate at $s=i$:
\[ \rho(\omega_1)|_{s=i}\ =\rho(\omega_3)|_{s=i}\
= i \left( \begin{array}{cc} 1& 1\\ 0 & 1
  \end{array}\right) , \ \ \ \ \ \  \rho(\omega_2)|_{s=i}\
= i \left( \begin{array}{cc} 1& 0\\ -1 & 1
  \end{array}\right)~.\]

Thus, we see that $\rho$ evaluated at $s=i$ is, up to powers of $i$,
equal to the
homology representation $h : M \ra \PSL_2(\Z)$. Another way to
state this is that   $\rho$
evaluated at $s=i$ followed by the projection $\GL_2(\C) \ra \PSL_2(\C)$
is the same as $h$ followed by the inclusion $\PSL_2(\Z) \ra
\PSL_2(\C)$.

Since $h$ determines the Nielsen-Thurston classification (Lemma~\ref{L1}),  we have
proved:

\begin{thm} The Nielsen-Thurston classification of a mapping class
  $\phi\in M=M(0,4)$ is determined, for every $n\geq 2$, by the
collection of quantum
  representations
$\rho_{n,k}(\phi)$ ($k\geq 1$). Moreover, in the Pseudo-Anosov case, these representations also determine the stretching factor.
\end{thm}

An immediate consequence is that for all $n$, we have
\[\bigcap_{k=1}^{\infty} \ker(\rho_{n,k}) = \ker \rho = \ker h = N~.\]
We can
strengthen this statement by considering the associated projective
representations, which we will now consider as homomorphisms $M\ra
\PSL_2(\C)$ and denote by $\rho_{n,k}^P$.
\begin{cor}\label{asf}
For all $n\geq 2$, we have that
\[\bigcap_{k=1}^{\infty} \ker(\rho^P_{n,k}) = N.\]
\end{cor}
This is the appropriate analog of the asymptotic faithfulness theorem
of the first author in \cite{A3} (see also Freedman, Walker, and
Wang \cite{FWW} for a proof in the skein-theoretical context).
\begin{proof} If $\phi$ is in the kernel of all $\rho^P_{n,k}$, then the
  off-diagonal entries of all matrices $\rho_{n,k}(\phi)$ are
  zero. It follows that the same is true for
  $\rho(\phi)$, hence also for $h(\phi)$, proving that
  $h(\phi)=1$.   \end{proof}

\section{The homology representation as a limit of quantum
  representations}
\label{sec4}

In the preceding section, we have exhibited the homology
representation $h$ as the evaluation at $s=i$ of the universal quantum
representation
$\rho$ with coefficients in $\Z[s,s^{-1}]$. Observe that the evaluation $s=i$ corresponds to $q=-1$. If we
consider the skein-theoretical construction of $\rho_{n,k}$
instead of the construction via conformal blocks, we can take for
$q$ not just  $ e^{2 \pi i/(k+n)}$ but any
root of unity of order $k+n$. In particular, we can choose a sequence
of roots of unity which converges to $-1$ as the level goes to
infinity,
and thereby exhibit the homology representation $h$ as a limit of the
quantum representations. This has interesting applications (see
Theorem~\ref{5.1} and Corollaries \ref{cor1}, \ref{cor12},  and \ref{cor2}).
For simplicity,
we restrict to the $SU(2)$-case, using the skein-theoretical approach to the
$SU(2)$-TQFT of \cite{BHMV2}.

The skein-theoretical analog of the representation $V_{2,k}$ is the
Kauffman bracket skein module of a $3$-ball with four marked points
colored $1$ on
the boundary. It can be
defined over the Laurent polynomial ring $\Z[A,A^{-1}]$. If one
then specializes $A$ to be a primitive $(4k+8)$-th root of unity, one
gets exactly the $SU(2)$-TQFT vector space at level $k$.
It has a standard basis
given by the tangles

\setlength{\unitlength}{1mm}
\centerline{
\begin{picture}(60,13)
\thicklines
\put(15,3){$h_0=$}
\spline(30,0)(35,5)(30,10)
\spline(40,0)(35,5)(40,10)
\end{picture}
\begin{picture}(60,13)
\thicklines
\put(15,3){$v_0=$}
\spline(30,0)(35,5)(40,0)
\spline(30,10)(35,5)(40,10)
\end{picture}}
\vskip 8pt

As usual, the lines in this figure stand for bands (or ribbons), and
the four boundary points must also be thought of as ``banded'' points
({\em i.e.,} small intervals). (This corresponds to the
 projective tangent
vectors used in the geometric construction.)

We represent a generator $\sigma_i$ of $\tM$ by the corresponding obvious half-twist of the
$3$-ball, corrected by small half-twists near the two banded
points permuted by $\sigma_i$ so as to preserve the band structure at
these points.  We can arrange the correction so that
$\sigma_1$ and $\sigma_3$ leave the tangle $h_0$ fixed, and  $\sigma_2$ leaves
the tangle $v_0$ fixed. This gives a representation of $\tM$
which we denote by $\tilde\rho^{(S)}$ (where the $S$ stands for {\em skein
theory}). It will be convenient to write our
matrices in the basis  $(A^{-1} h_0, A v_0)$. Then a Kauffman bracket calculation
gives the following matrices:

\[\tilde\rho^{(S)}(\sigma_1) =  \tilde\rho^{(S)}(\sigma_3) = \left( \begin{array}{cc} 1& -A^4\\ 0 & -A^4
  \end{array}\right)~, \ \ \  \tilde\rho^{(S)}(\sigma_2) = \left( \begin{array}{cc} -A^4& 0\\ -1& 1
  \end{array}\right).\]

We omit the details, as the calculation is similar to the one in \cite{M}, where the second author already computed the matrix for the
mapping class $\sigma_1^{-2}\sigma_2^{2}$ (and showed that this matrix has infinite order in the
$SU(2)$-TQFT representations at level $k$ except if $k=1,2,4,$ or
$8$).

This $\tilde\rho^{(S)}$ is only a projective representation of $M$, but if we define
\[\rho^{(S)}(\omega_i) = A^{-2} \tilde\rho^{(S)}(\sigma_i)~,\] we get a
linear representation of $M$. Moreover, $\rho^{(S)}$ evaluated
at $A^2=-i$ coincides with $\rho$ evaluated at $s=i$, where $\rho$ is
the universal quantum representation we considered in
the geometric context. Recall that this specialization is essentially
the homology representation $h$.

Let us now use this to exhibit $h$ as a limit of $SU(2)$-TQFT
representations. We need the following simple lemma.

\begin{lem} For every integer $r\geq 3$ there exists an integer
  $\ell=\ell(r)$ such that $\gcd(l,4r)=1$ and $\vert 2\ell - r\vert
 \leq 4$.
\end{lem}

\begin{proof} Take $\ell=(r+2)/2$ (respectively $(r+1)/2$, $(r+4)/2$,
 $(r-1)/2$) according to whether $r$ is congruent to $0$ (respectively
 $1,2,3$) modulo $4$.
\end{proof}

Let $A_k=e^{-2 \pi
  i\ell/(4k+8)}$ where $\ell=\ell(k+2)$. This is a primitive $(4k+8)$-th root of unity so by the
construction of \cite{BHMV2} the
specialization $A=A_k$ gives rise to a TQFT.

\begin{rem}\label{r5} {\em For this choice of root of unity (and big enough $k$), this TQFT is not unitary. In other words, the natural non-degenerate hermitian form on the TQFT vector spaces is (in general) indefinite.
Compare \cite[Remark 4.12]{BHMV2}.
}\end{rem}

Let
us denote by $\rho_k^{(S)}$ the corresponding representation
of $M$ as described above. Everything
has been done so that
$A_k^2$ converges to $-i$ as $k \ra \infty$. Thus for every mapping
  class $\phi\in M$, the limit of the matrices $\rho_k^{(S)}(\phi)$
  exists and is equal to a power of $i$ times $h(\phi)$. This proves the following

\begin{thm}\label{5.1} For every mapping class $\phi\in M=M(0,4)$, we have

\[\lim_{k\ra \infty} \vert \Tr \rho_k^{(S)}(\phi)\vert = \vert
\Tr
h(\phi)\vert ~.\]
\end{thm}

\begin{rem}{\em We may replace $\rho_k^{(S)}$ by the original
  projective representation $\tilde\rho_k^{(S)}$ in this theorem, as long as
  we agree that $\tilde\rho_k^{(S)}(\phi)$ means any matrix obtained by
  writing $\phi$ as a word in the generators and applying
  $\tilde\rho_k^{(S)}$ letter by letter to this word. This is because the
  only projective ambiguities of  $\tilde\rho_k^{(S)}$ are roots of
unity, so that
  one has $\vert \Tr  \tilde\rho_k^{(S)}(\phi)\vert = \vert \Tr
   \rho_k^{(S)}(\phi)\vert$.}\end{rem}

The following corollary generalizes the result
of \cite{M}.
\begin{cor}\label{cor1} Let $\phi\in M=M(0,4)$ be Pseudo-Anosov. Then no power
  of the
  matrix $\rho_k^{(S)}(\phi)$ is a multiple of the identity matrix,
  except possibly for finitely many values of $k$.
 \end{cor}
\begin{proof} Observe that the matrices $\rho_k^{(S)}(\omega_i)$ have
  determinant $-1$. Thus if a power of the matrix $\rho_k^{(S)}(\phi)$ is a
multiple of
the identity, then that multiple is $\pm 1$ or $\pm i$. Therefore a
further power of $\rho_k^{(S)}(\phi)$ is the identity, hence the eigenvalues
of $\rho_k^{(S)}(\phi)$ must be roots of unity and thus $\vert \Tr
  \rho_k^{(S)}(\phi)\vert\leq 2$. But for $k$ big enough this
  is impossible by our theorem, since $\phi$ Pseudo-Anosov
  implies $\vert
\Tr
h(\phi)\vert \geq 3$.
 \end{proof}
\begin{cor}\label{cor12} A mapping class $\phi\in M=M(0,4)$ is
  Pseudo-Anosov if and only if its $SU(2)$ TQFT representation matrix
  $\rho_k^{(S)}(\phi)$ has infinite order for large enough level $k$.
\end{cor}
\begin{proof} It suffices to show that reducible mapping classes are
  represented by matrices of finite order. This is easily checked, as
 every reducible mapping
  class is a conjugate of a power of $\omega_1$ times an element of the translation
 subgroup $N\cong
  (\Z/2\Z)^2$ (this follows from the description in Corollary 3.3).
\end{proof}
\begin{rem} {\em Corollaries \ref{cor1} and \ref{cor12} are algebraic
statements, and
  therefore true also for the level $k$ TQFT where $A$ is an arbitrary
  primitive $(4k+8)$-th root of unity.} \end{rem}

In a similar vein, we have the following corollary (here we
assume that $A$ is the root $A_k$ as defined above).
\begin{cor}\label{cor2} Let $\phi\in M=M(0,4)$ be Pseudo-Anosov with
  stretching factor $\lambda(\phi)$. Then for $k$ big enough
$\rho_k^{(S)}(\phi)$ has a unique eigenvalue $\lambda_k$ such that
$\vert\lambda_k\vert >1$, and \[\lim_{k\ra \infty} \vert\lambda_k\vert
=  \lambda(\phi)~.\]
\end{cor}
 \begin{proof} The eigenvalues of $\rho_k^{(S)}(\phi)$ converge to a
   power of $i$ times
   those of $h(\phi)$. Since $\vert \Tr h(\phi) \vert =
   \lambda(\phi) + \lambda(\phi)^{-1}$, the result
   follows. \end{proof}

Note that we cannot obtain such a result from a unitary
TQFT (such as the one we get if we put $A=e^{2 \pi i/(4k+8)}$)
 because there the eigenvalues of the matrices
 $\rho_k^{(S)}(\phi)$ would have unit norm.

\begin{rem} {\em A similar story exists for the mapping class group
   $M(1,1)$ \cite{M2}.}\end{rem}

\end{document}